\documentclass{article}       
\usepackage{fullpage}
\usepackage{changepage}  

\usepackage{amsthm} 
\usepackage{enumitem} 
\usepackage{hyperref}
\usepackage{amssymb} 
\usepackage{amsmath} 
\usepackage{bm} 
\usepackage{mathtools}   
\usepackage{booktabs} 

\usepackage{tikz}
\usetikzlibrary{positioning,backgrounds,fit,chains,scopes,shapes,arrows}
\usepackage{pgfplotstable} 

\newcommand{\onetwo}{\textstyle{\frac{1}{2}}}
\DeclareMathOperator*{\argmin}{arg\, min}

\DeclareMathOperator{\prox}{prox}

\newcommand{\mX}{{U}}
\newcommand{\mY}{{V}}
\newcommand{\ind}{\imath_{\mathbb{R}_+}}

\newcommand{\tv}{\mathrm{tv}}
\newcommand{\gt}{u_{\text{gt}}}
\newcommand{\diver}{\mathrm{div\, }}
\newcommand{\reg}{\mathcal{R}} 
\newcommand{\R}{\mathbb{R}}

\newtheorem{prop}{Proposition} 
\theoremstyle{remark}

\theoremstyle{definition} 
\newtheorem{alg}{Algorithm}

\title{\raggedright \textbf{FISTA-Condat-Vu: Automatic Differentiation for Hyperparameter Learning in Variational Models}  \\ \vspace{20pt} \normalsize Patricio Guerrero$^{*}$, Simon Bellens, Wim Dewulf \\
Department of Mechanical Engineering, KU Leuven, Celestijnenlaan 300, 3001 Leuven, Belgium and M\&A Corelab, Flanders Make, KU Leuven  \\
\medskip
$^*$Corresponding author: Patricio Guerrero, patricio.guerrero@kuleuven.be
\vspace{-30pt}}
\author{}
\date{}

\begin{document}
\maketitle

\begin{adjustwidth}{30pt}{30pt}
\noindent\textbf{\small Abstract.} 
{\small Motivated by industrial computed tomography, we propose a memory efficient strategy to estimate the regularization hyperparameter of a non-smooth variational model. The approach is based on a combination of FISTA and Condat-Vu algorithms exploiting the convergence rate of the former and the low per-iteration complexity of the latter. The estimation is cast as a bilevel learning problem where a first-order method is obtained via reduced-memory automatic differentiation to compute the derivatives. The method is validated with experimental industrial tomographic data with the numerical implementation available. \\
\textbf{keywords}: computed tomography, few-view, regularization parameter, automatic differentiation, bilevel learning  
}
\end{adjustwidth}

\section{Introduction}

We consider ill-posed imaging inverse problems of the form
\begin{equation}
   v = Au + \text{noise}
\label{eq.inv_prob}
\end{equation}
where the \emph{forward operator} $A\colon \mX \to \mY$ is a linear operator between real Hilbert spaces $\mX$ and $\mY$; $u\in \mX$ is the image to recover and $v\in \mY$ is the data modelled by $A$. 

Problem (\ref{eq.inv_prob}) is ill-posed in the sense that the solution does not depend continuously on the data, if it exists. A common strategy is to look for regularized solutions by solving the optimization problem
\begin{equation}
    \min_{u \in \mX} f(u)+\lambda  \reg(u), \quad \lambda > 0,
    \label{eq.optimization}
\end{equation}
where $f$ is a \emph{data fidelity} term that measures some distance or error between $Au$ and the data $v$; $\reg$ is a suitable \emph{regularizer}. Our goal is to estimate, via bilevel learning, $\lambda$, a regularization parameter that balances the influence of each term in (\ref{eq.optimization}).

We will make the following assumptions on the involved functions. Let $\bar\R = \R \cup \{+\infty\}$,
\begin{enumerate}[label={($a$}{{\arabic*}}),leftmargin=3\parindent]
    \item $f\colon \mX \to \R$ is a convex and differentiable function with a $\beta-$Lipschitz continuous gradient for some $\beta>0$.
    \item $\reg\colon \mX \to \bar\R$, is a convex, proper, lower semi-continuous (lsc) function, allowed to be non-smooth. 
    \item Problem (\ref{eq.optimization}) admits at least one solution for any $\lambda > 0$.
\end{enumerate}

\paragraph{Notations} $\mX$ and $\mY$ will be endowed with an inner product $ \langle \cdot, \cdot \rangle$ and induced norm $\| \cdot \| = \langle \cdot, \cdot \rangle ^{1/2}$. The gradient operator of a 3D image is the linear operator $\nabla \colon \mX \to W = \mX\times\mX\times\mX$, $\nabla = (\nabla_x, \nabla_y, \nabla_z)$ of partial derivatives of the image. The adjoint of a bounded linear operator $B\colon \mX \to W$ is the operator $B^* \colon W \to X$ such that 
	$\langle Bu , w\rangle = \langle u,  B^* w\rangle$ for all $u\in \mX$ and $w\in W$.

Relevant examples are $\nabla^* = -\diver$ (the negative divergence operator) \cite{chammbolle_pock2011} and for a tomographic projection operator $A$, $A^*$ is the \emph{backprojection} operator \cite{guerrero2023}. The (isotropic) \emph{total variation} (TV) functional of a 3D image $u$ will be denoted by $\tv(u)$, which, after discretization is defined as the $\ell_{2,1}$ norm of the gradient of the image, i.e., $\tv(u) = \| \nabla u \|_{2,1} = \| |\nabla u|_2 \|_1$ where $|\cdot|_2$ is the component-wise Euclidean norm, see \cite{chambolle_pock2016} for details. The \emph{indicator} function of a set $\Omega$, denoted by $\imath_\Omega$ is defined as $\imath_\Omega(u) = 0$ if $u\in \Omega$ and $\imath_\Omega(u)=+\infty$ if $u\notin\Omega$. The \emph{Heaviside} step function $H$ is defined as $H(\lambda)=0$ if $\lambda<0$ and $H(\lambda)=1$ if $\lambda\geq0$. The \emph{ramp} function is defined as $\mathrm{ramp}(\lambda) = 0$ if $\lambda < 0$ and $\mathrm{ramp}(\lambda) =  \lambda$ if $\lambda\geq0$.

\paragraph{Problem setting} We seek to find $\lambda$ in (\ref{eq.optimization}) that minimizes a convex and differentiable loss function $L$. We assume that we dispose of \emph{at least} one ground-truth image $\gt$ related to some corrupted data $v$ to perform the learning. Then, we are interested in $\lambda^*$ that solves the bilevel optimization problem
\begin{equation}\label{eq.bilevel}
\begin{array}{ll}
  \lambda^* = &\argmin\limits_{\lambda >0} \left\{ L(\lambda) \coloneq \onetwo \| u_\lambda - \gt \|^2 \right\} 
  \\
  & \text{subject to} \ u_\lambda \in \argmin\limits_{u\in\mX}  f(u)+\lambda \reg(u),
\end{array}
\end{equation}

\paragraph{Motivating application} The strategies developed are motivated by an industrial computed tomography reconstruction problem with incomplete data, detailed in Section \ref{sec.numerics}. A common approach is to adapt Problem (\ref{eq.optimization}) with a squared $L^2$ norm data fidelity and TV regularization enforcing piecewise constant solutions. In addition, we will combine TV with the indicator function of the non-negative orthant, enforcing  non-negative solutions, that is, we will consider in our application 
\begin{equation}
  f= \onetwo\|A \cdot - v \|^2, \quad \reg = \tv + \imath_{\mathbb{R}_+}.
  \label{eq.functions}
\end{equation}

 However, the methodology presented here is also valid for arbitrary $f$ and $\reg$ verifying assumptions $(a1)-(a3)$ or an arbitrary convex and differentiable loss $L$. In fact, it will be useful for the derivation of the algorithms to \emph{split} $\reg$ as the sum and composition $\reg = g \circ L + h$ with $L$ a bounded linear operator and focus instead on the problem
\begin{equation}
    \min_{u \in \mX} f(u)+\lambda \left( g(Lu) + h(u)\right), \quad \lambda > 0.
    \label{eq.optimizationgenearal}
\end{equation}

Finally, in (\ref{eq.functions}) we have 
$
      g = \|\cdot \|_{2,1}, \  L=\nabla, \  h = \ind.
$

\paragraph{Related work} Results on hyperparameter estimation for Problem (\ref{eq.optimization}) include the discrepancy principle \cite{morozov, per},  L-curve \cite{Lcurve}, GCV \cite{gcv} or SURE \cite{sure}. These methods are either empirical, need the knowledge of the noise level or rely on a brute force grid-search strategy. Bilevel learning emerged as an alternative in the imaging community, e.g. \cite{bilevelreyes, kunisch2013bilevel, optimalparameter}. These approaches usually require a smoothing of the regularization term and are solved with implicit differentiation. To our knowledge, this is the first work related to memory usage in automatic differentiation (AD) to solve such bilevel problems without any smoothing of the involved functions.

\paragraph{Contributions}

The main contribution of this work is a computationally efficient strategy to solve bilevel Problem (\ref{eq.bilevel}) for 3D imaging problems based on AD. 
The solution of such a problem is the regularization hyperparameter of a non-smooth variational model. In particular, we show the equivalence of Condat-Vu (CV) and gradient projection (GP) algorithms for a denoising problem, with a key difference in their memory footprint when AD is performed (with respect to the regularization parameter). Then, we propose aCV (assisted CV) where we reduce even more the memory storage related to GP. Then, we address the parameter learning problem with a combined FISTA-aCV exploiting the optimal convergence rate of FISTA and the low per-iteration memory footprint of aCV. This results in an accelerated first-order algorithm that is tested on an industrial case with the python implementation available at \texttt{github.com/patoguerrero/fista\_acv}. 

\section{Image reconstruction}\label{sec.lower}

In this preliminary section, we will explore solutions of (\ref{eq.optimizationgenearal}) for a fixed $\lambda$ via proximal methods, which are algorithms based on the proximity operator of the functions to minimize. The combination of FISTA and Condat-Vu will be introduced here.

 We recall the definition of the \emph{proximity operator} \cite{bauschke_book2017} of a convex, proper, lsc function $g$, defined for any $u\in\mX$ by $\prox_g (u) = \argmin\limits_{w\in\mX} \Big\{ \onetwo \| u -  w \|^2  + g(w) \Big\}$, which exist and is unique for any $g$ \cite{condat_review2023,chambolle_pock2016}. Also, we recall the \emph{Legendre–Fenchel conjugate} associated to any function $g\colon \mX \to \bar\R$ defined for $u\in\mX$ by
    $g^*(u) = \sup\limits_{w\in\mX} \langle u,w \rangle  -  g(w)$,
which is a convex and lsc function. 



\subsection{FISTA and Condat-Vu}

In \cite{mfista} an accelerated proximal gradient descent algorithm referred to FISTA is proposed and provides an optimal convergence rate of  $O(\beta/k^2)$ for convex functions \cite{CVacceleration}. The resulting iteration adapted to Problem (\ref{eq.optimization}) is, for all $k\geq0$ 
\begin{equation*}
\text{FISTA :} \quad \left\{
\begin{array}{ll}
  & u_{k+1} = \prox_{\gamma\lambda \reg} \left(\hat u_k - \gamma \nabla f(\hat u_k) \right),  \\
  & t_{k+1} = \onetwo \left( 1 + \sqrt{1+4t_k^2}  \right), \\
  & \hat u_{k+1} = u_{k+1} + \dfrac{t_k - 1}{t_{k+1}} (u_{k+1} - u_{k}), 
\end{array} 
\right.
\end{equation*}
with $t_0 = 1$, $u_0, \hat u_0 \in \mX$ and a step size $\gamma\in(0,1 / \beta]$. 

Note that the function $\prox_{\gamma\lambda\reg}$ needs to be computed at each iteration. Therefore, this method is attractive especially when $ \prox_{\reg} $ has either a closed-form expression or it can be computed efficiently. 

The Condat-Vu (CV) algorithm, was proposed in \cite{condat2013, vu2013} to solve problems in the form of (\ref{eq.optimizationgenearal}) with the iteration
\begin{equation*}
\text{CV :} \quad \left\{
\begin{array}{ll}
  & u_{k+1} = \prox_{\gamma\lambda h} \left(u_k - \gamma L^* w_k - \gamma \nabla f(u_k) \right),  \\
  & w_{k+1} = \prox_{\sigma(\lambda g)^*} \left( w_k + \sigma L u_{k+1} \right), 
\end{array}
\right.
\end{equation*}
with $u_0 \in \mX$, $w_0 \in W$. The step sizes $\gamma>0$ and $\sigma>0$ are required in \cite{condat_review2023} to verify $\gamma\sigma\| L \|^2 <  1$ and $\gamma < 1/ \beta$ and then the iteration converges (weakly, see \cite{condat_review2023}) to a solution. 



\subsection{Preliminary application: the proximity operator of $\mathrm{tv} + \ind$} This example is motivated by the need of computing $\prox_{\lambda\reg}$ in FISTA. Given $v\in\mX$ and a fixed $\lambda > 0$, computing $\prox_{\lambda\reg}(v)$ with $\reg = \mathrm{tv} + \ind$ results in solving
\begin{equation}\label{eq.denoising}
    \argmin\limits_{v\in\mX} \left\{ \onetwo \| u - v \|^2  + \lambda \tv(u) + \ind(u) \right\},  
\end{equation}
which we can immediately identify as a standard TV \emph{denoising} problem for a noisy image $v$ with a non-negativity constraint. The algorithms \emph{gradient projection} (GP) and \emph{fast gradient projection} (FGP) were proposed in the original work \cite{mfista} for that goal, GP iterates as
\begin{equation*}
\text{GP :} \quad \left\{
\begin{array}{ll}
  & u_{k} = \Pi_{\R_+} \left( v + \lambda \diver w_k \right),  \\
  & w_{k+1} = \Pi_{\{ \| \cdot \|_{2, \infty} < 1 \}} \left( w_k + 1/(\lambda\|\nabla\|^2) \nabla u_{k} \right), 
\end{array}
\right.
\end{equation*}
whereas FGP accelerates GP on $w_k$. Here $\Pi_{\R_+}$ is the (component-wise) orthogonal projection onto $\R_+$ and $\Pi_{\{ \| \cdot \|_{2, \infty} < 1 \}}$ is the projection onto 2-balls with radius $1$, computed e.g. as in \cite[eq. (4.23)]{chambolle_pock2016}, expressed in (\ref{eq.projection}) below.

On the other hand, note that adapting CV for the denoising problem (\ref{eq.denoising}) results in an iteration where we can observe that $\lambda$ acts only on $\prox_{(\lambda g)^*}$ and not in $\prox_{\lambda h}$ as $\lambda\ind = \ind$ while in GP it acts on both updates.

\subsection{Equivalence between GP and CV for denoising} 

\begin{prop}
    GP is a dual instance of CV applied on the TV denoising Problem (\ref{eq.denoising})
\end{prop}
This is a trivial exercise if we supply CV with the operators $A=I$, $L=\nabla$, $h=\ind$, $f = (1/2) \|\cdot - v \|^2$, $g = \| \cdot \|_{2,1}$ (then $g\circ L = \tv$) and the step sizes
\begin{equation}
    \gamma = \lambda, \quad \sigma=1/(\lambda\|\nabla\|)^2.
\label{eq.stepGP}
\end{equation}

These (and only these) step sizes allow to \emph{remove} the previous $u_k$ dependency in the update $u_{k+1}$ and therefore GP becomes a dual algorithm from the primal-dual CV. 

Note that to be able to write $\gamma = \lambda$, $1/\lambda$ needs to be the Lipschitz constant of the gradient of $f$, therefore problem (\ref{eq.denoising}) needs to be rewritten accordingly, particularly, by letting $f_\lambda = (1/\lambda) f$ and then the problem GP is solving is actually 
\begin{equation}
\min_{u \in \mX}  f_\lambda(u)+  g(Lu) + h(u), \quad \lambda > 0. 
\label{eq.gpmin}
\end{equation}

Importantly, even if GP and CV are equivalent to solve (\ref{eq.denoising}), in GP $\lambda$ is included in $f$ and not in $\reg$ as is in the version of CV presented here, this makes them equivalent algorithms \emph{with different step sizes}. Observe that both step sizes in (\ref{eq.stepGP}) depend on $\lambda$.

Finally, CV applied to  Problem (\ref{eq.gpmin}) requires the step sizes to actually verify $\gamma < \lambda$ and $\sigma<1/(\lambda\|\nabla\|^2)$ to guarantee convergence in the sense detailed in \cite{condat_review2023}. The algorithm presented there is a \emph{relaxed} algorithm, which means that there is another auxiliary variable helping to improve the convergence rate. Although this is possible here, we do not employ relaxation mainly because this will increase the memory footprint of CV. Therefore, we can tune the step sizes in CV as (\ref{eq.stepGP}) which makes GP and CV equivalent. The convergence of the iteration is still assured now in the sense of \cite{mfista}. For the same reason, CV applied now to (\ref{eq.denoising}) can be tuned with 
\begin{equation}
    \gamma =1, \quad \sigma = 1 / \|\nabla\|^2,
    \label{eq.stepCV}
\end{equation}
where we can appreciate the non-dependency on $\lambda$. This fact will be exploited in the next section.

\section{Main result: hyperparameter learning}

Recall that we want to explore the estimation of $\lambda$ in (\ref{eq.optimization}) by solving the bilevel problem (\ref{eq.bilevel}). This will be done via (accelerated) gradient descend where the challenge is to compute the derivative $\nabla_\lambda L$ of $L$ with respect to $\lambda$. We will use AD in \emph{reverse mode} (rAD) as the main tool. rAD has the drawback of a large memory requirement making the tool unfeasible for high-resolution 3D imaging problems. Therefore, we propose a memory efficient algorithm to compute $\nabla_\lambda L$. 

\subsection{Derivative of the loss function} 

For a given algorithm that solves (\ref{eq.optimizationgenearal}), in rAD we need to store a certain number of variables that increases linearly with the iterations. Then, we propose to use FISTA to exploit its optimal convergence rate. However, FISTA needs the solving of (\ref{eq.denoising}). We will combine FISTA with CV, referred to the FISTA-CV algorithm, aiming at reducing the memory usage. Therefore, it is crucial that CV stores variables per iteration in rAD such that the problem is feasible in practice, e.g., with a $\sim$24 GB GPU card. It is worth to recall that even if GP and CV are equivalent, they are written differently and we will show below how this makes CV more memory efficient than GP. Even better, CV allows an easy improvement if we are able to inform the algorithm with the closed-form expression of the derivate of $\prox_{(\lambda g)^*}(w)$ with respect to $\lambda$. We will see that in many cases with non-smooth $g$, e.g. TV (where $g = \|\cdot\|_{2,1}$), this is indeed possible.   
\begin{figure}
\centering\footnotesize
\begin{tikzpicture}[%
    start chain=going right,    
    node distance=30pt and 20pt,
    every join/.style={norm},   
    ]
\tikzset{
  base/.style={draw, rounded corners, on chain, on grid, minimum height=5ex, align=center},
  nada/.style={on chain, on grid, minimum height=5ex, align=center},
  azul/.style={base, rounded corners, fill=blue!20},
  norm/.style={->, draw},
}
\node [base] (wk) {$w_k$};
\node [base, join, right =12ex of wk] (div) {$\lambda\diver w_k + v$};
\node [base] (uk) {$u_k$};
\node [base, join, right =16.5ex of uk] (wk2) {$w_k + 1/(\lambda\|\nabla\|^2)\nabla u_k$};
\node [base] (wk3) {$w_{k+1}$};
\node [azul, above=of wk] (lbd) {$\lambda$};
\node [nada, below=4ex of div] (a) {$g_1$};
\node [nada, below=4ex of wk2] (b) {$g_2$};

\draw [->] (div) -- node[above] {$\Pi_{\R_+}$} (uk);
\draw [->] (wk2) -- node[above] {$\Pi_1$} (wk3);
\draw [->] [out=0, in=90] (lbd) to (div);  
\draw [->] (lbd.east) [out=0, in=155] to (wk2.north);
\draw [->] (wk.north east) [out=25, in=155] to (wk2.north west);
\end{tikzpicture}
\caption{Computational graph of a GP iteration: the map $w_k \mapsto w_{k+1}$.}
\label{fig.graphGP}
\end{figure}
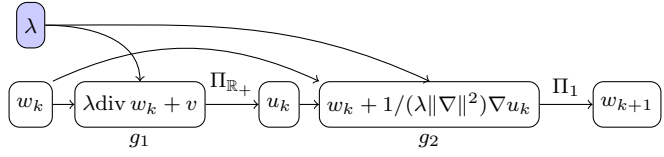
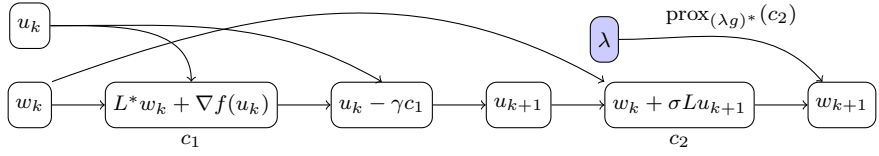
\begin{figure}
\centering\footnotesize
\begin{tikzpicture}[%
    start chain=going right,    
    node distance=30pt and 20pt,
    every join/.style={norm},   
    ]
\tikzset{
  base/.style={draw, rounded corners, on chain, on grid, minimum height=5ex, align=center},
  nada/.style={on chain, on grid, minimum height=5ex, align=center},
  azul/.style={base, rounded corners, fill=blue!20},
  norm/.style={->, draw},
}
\node [base] (wk) {$w_k$};
\node [base, join] (div) {$ L^* w_k + \nabla f(u_k)$};
\node [base, join] (uk2) {$u_k - \gamma c_1$};
\node [base, join] (uk3) {$u_{k+1}$};
\node [base, join] (wk1) {$w_k + \sigma L u_{k+1}$};
\node [base, join] (wk2) {$w_{k+1}$};
\node [azul, above=16pt of wk1.west] (lbd) {$\lambda$};
\node [base, above=of wk] (uk) {$u_k$};
\node [nada, below=4ex of div] (a) {$c_1$};
\node [nada, below=4ex of wk1] (b) {$c_2$};

\draw [->] [out=0, in=130] (lbd) to node[above] {$\prox_{(\lambda g)^*}(c_2)$} (wk2);  
\draw [->] [out=0, in=90] (uk) to (div); 
\draw [->] [out=0, in=150] (uk) to (uk2.north);
\draw [->] (wk.north east) [out=20, in=150] to (wk1.north west);
\end{tikzpicture}
\caption{Computational graph of a CV iteration: the map $(u_k,w_k) \mapsto (u_{k+1},w_{k+1})$. Note that the constants $\{\gamma, \sigma\}$ do not depend on $\lambda$.}
\label{fig.graphCV}
\end{figure}

\paragraph{Sensitivity of $\lambda$ in GP and CV} Let us abbreviate the projection operator onto 2-balls used in GP/CV as $  \Pi_\lambda \coloneq \Pi_{\{ \| \cdot \|_{2, \infty} < \lambda \}} $. We add the following assumption to $(a1)-(a3)$:
\begin{enumerate}[label={($a$}{{\arabic*}}),leftmargin=3\parindent]
\setcounter{enumi}{3}
\item $h$ in (\ref{eq.optimizationgenearal}) is such that $\lambda h = h$ for any $\lambda>0$.
\end{enumerate}

If $\reg = \tv + \ind$, we have $h=\ind$ which verifies $(a4)$ as an indicator function of any set does. Diagrams \ref{fig.graphGP} and \ref{fig.graphCV} show the (summarized) computational graph of a single iteration of GP and CV under $(a4)$ respectively, that is, the maps $w_k \mapsto w_{k+1}$ and $(u_k,w_k) \mapsto (u_{k+1},w_{k+1})$ where we aim at the sensitivity of the output variables with respect to $\lambda$. When analyzing FISTA with GP as originally proposed in \cite{mfista}, we observe that $\lambda$ acts on both step sizes to obtain intermediate variables $g_1$ and $g_2$. Therefore, from Graph \ref{fig.graphGP} and using differential calculus, the sensitivity of $\lambda$ on $w_{k+1}$ is computed by
\begin{equation}
    \frac{\partial w_{k+1}}{\partial \lambda} = 
    \frac{\partial g_1}{\partial \lambda}
    \frac{\partial u_{k}}{\partial g_1}
    \frac{\partial g_2}{\partial u_k}
    \frac{\partial w_{k+1}}{\partial g_2} + 
    \frac{\partial g_2}{\partial \lambda}
    \frac{\partial w_{k+1}}{\partial g_2},
\label{eq.gpsensitivity}
\end{equation}
where at each partial derivative step, the storage in RAM memory of certain variables is needed in rAD. 

On the other hand, in CV under $(a4)$, $\lambda$ only intervenes in the last projection step. As a primal-dual algorithm we should be interested in the sensitivities $\dfrac{\partial u_{k+1}}{\partial \lambda}$ and $\dfrac{\partial w_{k+1}}{\partial \lambda}$ but as $u_{k+1}$ does not depend on $\lambda$, the remaining sensitivity comes from the following result. 
\begin{prop}
Under assumptions $(a1)-(a4)$, the sensitivity of $\lambda$ on a single CV iteration is  
\begin{equation}
    \frac{\partial w_{k+1}}{\partial \lambda} = 
    \frac{\partial}{\partial \lambda} \prox_{(\lambda g)^*}(c_2), 
    \label{eq.cvsensitivity}
\end{equation}    
provided such derivative exists for almost all $\lambda>0$. In particular, for $g=\| \cdot \|_{2,1}$ we have $\prox_{(\lambda g)^*} = \Pi_\lambda$, then 
\begin{equation}
\frac{\partial}{\partial \lambda} \Pi_\lambda(w) = \dfrac{w |w |_2 }{\lambda^2 \max^2\left\{1, \lambda^{-1}|w|_2\right\}} H(\lambda^{-1}|w|_2)),
\label{eq.diffprojection}
\end{equation}
with all operations to be understood component-wise.
\end{prop}

Expression (\ref{eq.cvsensitivity}) comes from Graph (\ref{fig.graphCV}) and differential calculus while expression (\ref{eq.diffprojection}) is obtained after writing the definition of $\Pi_\lambda$ \cite[eq. (4.23)]{chambolle_pock2016} which is
\begin{equation}
    \Pi_\lambda(w) = \dfrac{w}{\max \left\{1, \lambda^{-1} | w |_2 \right\}},
\label{eq.projection}
\end{equation}
where all operations are component-wise. Then, the map $\lambda\mapsto\max(1,\lambda)$ can be written with the ramp function as 
$\max(1,\lambda) = \mathrm{ramp}(\lambda-1)+1$, whose derivative is the Heaviside step function $H$ and (\ref{eq.diffprojection}) follows. We see that $\mathrm{ramp}$ is not differentiable at $0$ but then we can define $H(0)=0$ as the subgradient of $\mathrm{ramp}$ at $0$ is the interval $[0,1]$. 

Thus, computing the sensitivity of $\lambda$ within a single iteration of CV depends only on its sensitivity on $\Pi_\lambda (c_2)$ and it will be more memory efficient than (\ref{eq.gpsensitivity}) only if such a sensitivity requires fewer variables to be stored than those required in (\ref{eq.gpsensitivity}). This is trivial (and will be tested numerically) as (\ref{eq.gpsensitivity}) is decomposed into more elementary functions than (\ref{eq.projection}).

\paragraph{assisted Condat-Vu and the FISTA-aCV algorithm}

$\lambda \mapsto \Pi_\lambda(c_2)$ is not an elementary function and using AD to compute its derivative will still require to store a certain number of intermediate variables. Instead, we propose aCV, where we \emph{assist} CV and \emph{inform} the algorithm with its derivative. Under $(a4)$ and with $g = \|\cdot \|_{2,1}$, we can use (\ref{eq.diffprojection}) but the approach is valid for any $g$ having a closed-form expression of $\frac{\partial}{\partial \lambda} \prox_{(\lambda g)^*}$, for instance, $g = \|\cdot \|_{1,1}$ (anisotropic TV). Also, we could just take another $L$ in (\ref{eq.optimizationgenearal}), e.g., a shearlet transform \cite{tommy} or simply the identity operator where $g$ reduces to $\|\cdot \|_{1}$.

Providing the closed-form expression of the sensitivity (\ref{eq.diffprojection}) to a rAD chain results in two incompatible benefits: reducing the computing time or the memory footprint. The former is achieved by storing auxiliary variables in the forward pass also needed in the gradient computation, e.g., $\max\{1,\lambda^{-1}|w|_2\}$ or $|w|_2$, at the expense of memory usage. The latter is achieved by only storing $w$ in the forward pass within a single iteration of CV, to then compute the gradient using (\ref{eq.diffprojection}). This is the approach we used in our numerical implementation. Further improvement could be achieved by using compression techniques on $w$, or simply by storing $w$ with a lower bit depth.  Finally, FISTA combined with aCV to compute $\nabla_\lambda L$, here referred to FISTA-aCV, can be executed in any AD engine that allows an expansion of its elementary functions. 

\subsection{The resulting first-order algorithm for $\lambda$}

We are ready to present a first-order algorithm to solve Problem (\ref{eq.bilevel}) as we have the derivative $\nabla_\lambda L$ of the loss $L$ available. We will use a Nesterov-accelerated \cite{nesterov1983, chambolle_pock2016} gradient descent (NGD) algorithm where the step sizes $\gamma_k$ are selected with a line search approach following a sufficient decrease and backtracking strategy until the \emph{Armijo} condition is verified \cite[Algorithm~3.1]{nocedal}. Namely, if $\gamma_k$ does not verify the Armijo condition, we decrease $\gamma_k$ by a contraction factor of $1/2$. 

Paired with the Armijo step size selection and FISTA-aCV to compute the derivatives, the following NGD algorithm is a \emph{parameter-free} algorithm that summarizes our approach.

\vspace{10pt}
{\hrule width 0.48\textwidth}
\begin{alg}
\small
NGD for Problem (\ref{eq.bilevel}) \\
Given a pair $(v, \gt)$ of some data and a ground-truth reconstruction, $t_0=1$ and $\lambda_0,\hat\lambda_0>0$   
\begin{enumerate}[label=\footnotesize\arabic*), itemsep=0pt] 
    \item until convergence, for $k=0,1,\dots,$ do
    \begin{enumerate}[label=\footnotesize\arabic*)]\vspace{0pt}
    \item compute $d_k = \nabla_\lambda L({\hat\lambda_k}) $ with FISTA-aCV 
    \item $\lambda_{k+1} = \hat\lambda_{k} - \gamma_k d_k $, with $\gamma_k > 0$ from the Armijo rule 
    \item $t_{k+1} = \onetwo \left( 1 + \sqrt{1+4t_k^2}  \right)$ 
    \item $\hat \lambda_{k+1} = \lambda_{k+1} + \dfrac{t_k - 1}{t_{k+1}} (\lambda_{k+1} - \lambda_{k})$
    \end{enumerate}
    \item return $\lambda^* = \lambda_{k+1}$
\end{enumerate}
{\hrule width 0.48\textwidth}
\label{alg.agd}
\end{alg}

\section{Application: Industrial computed tomography}\label{sec.numerics}

\begin{table}[t]
\centering \small
\caption{Stored RAM memory in one GP/CV/aCV iteration for problem (\ref{eq.denoising}) related to an image of $1400^2\times2$ voxels.} 
    \begin{tabular}{llll}
    \hline\addlinespace[2pt]   
      & GP & CV & aCV \\   
       stored memory (MB.) & 153 & 82 & 49  \\ 
   \addlinespace[2pt]\hline\addlinespace[-4pt]      
    \end{tabular}
    \label{tab.memoryCV}
\end{table}

The \emph{cone-beam} transform models the scanning process in industrial CT, then it will be the considered forward operator $A$ in (\ref{eq.inv_prob}). See \cite{guerrero2024} for a detailed definition. Our experimental case is \emph{few-view} CT, i.e., when the number of acquired projections is drastically smaller than the number needed to satisfy the Nyquist criterium \cite{kak}. Then, a metal Ti-6Al-4V additively manufactured object was scanned with a Nikon XT-H-225-ST taking 3142 and 60 projections of $2000^2$ pixels. 

We implemented FISTA-GP, -CV and -aCV to compute $\nabla_\lambda L(\lambda)$ using the AD library \texttt{torch.autograd} of pytorch in python and monitored the RAM memory usage by the GPU within each iteration of GP, CV and aCV. The outputs are displayed in Table \ref{tab.memoryCV}. As expected we observe that CV reduces GP memory storage by 46$\%$ while aCV does it by 68$\%$.

The high-resolution scan was reconstructed with the standard FDK algorithm \cite{natteredwubbeling} and taken as our ground-truth reconstruction. Then, the few-view scan was reconstructed with FISTA-CV and Algorithm \ref{alg.agd} for the parameter learning. The only pre-processing task performed was the alignment of projections accounting for a detector horizontal shift and in-plane tilt, which was done with \texttt{alignCT} \cite{guerrero2024}.



While running FISTA-CV and -aCV, we observed convergence after 140 FISTA iterations with 5 inner CV iterations. Algorithm \ref{alg.agd} in turns showed convergence after 30 NGD iterations. The learning was performed with a volume of $1400^2\times2$ voxels. Under this setting, the learning was only possible with aCV and not with CV/GP using a 24 GB GPU card. The step sizes were selected as $\gamma = 1 /\beta $ for FISTA with $\beta$ computed with the power method \cite{sidky2012convex} and with (\ref{eq.stepCV}) for CV having $\| \nabla \| = \sqrt{18}$ (related to 3D images). The forward and backprojection operators $\{A, A^*\}$ were performed with \texttt{tomosipo} \cite{tomosipo} and the rest of the algebraic operations with pytorch, both in python.  

\begin{figure}[t]
    \centering
    \includegraphics[width=0.25\linewidth]{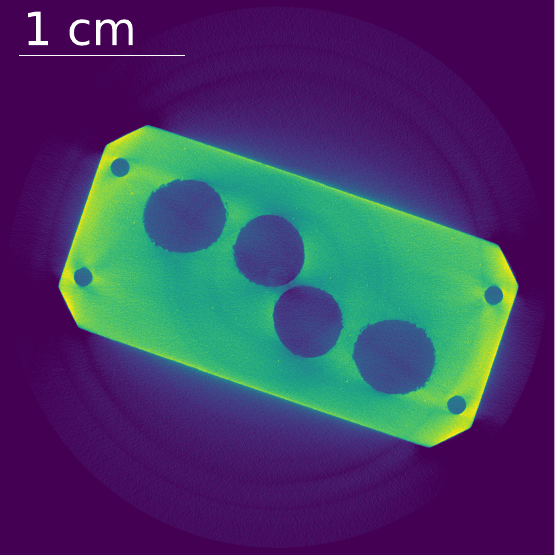}\hspace{-3pt}
    \includegraphics[width=0.25\linewidth]{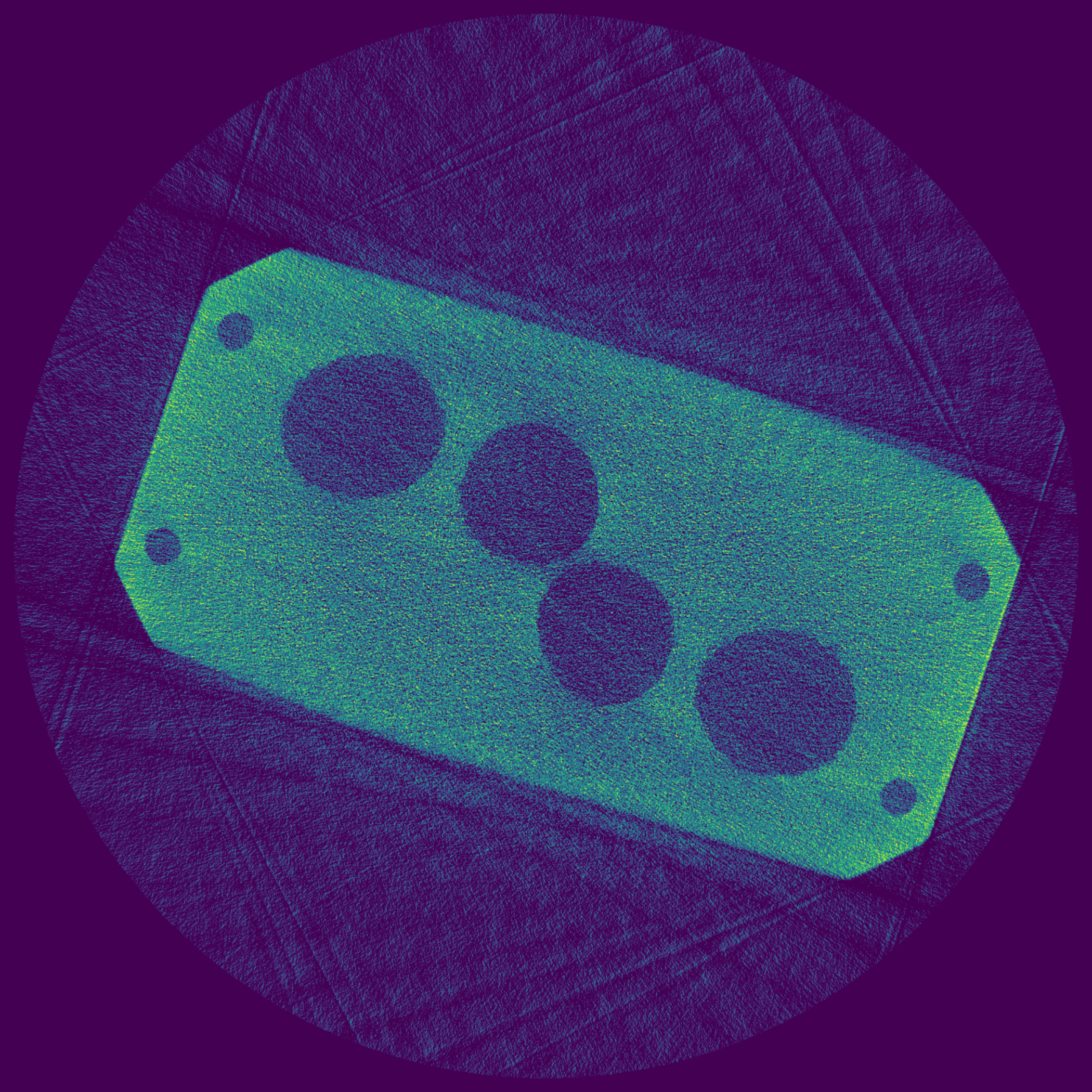}\hspace{-3pt}
    \includegraphics[width=0.25\linewidth]{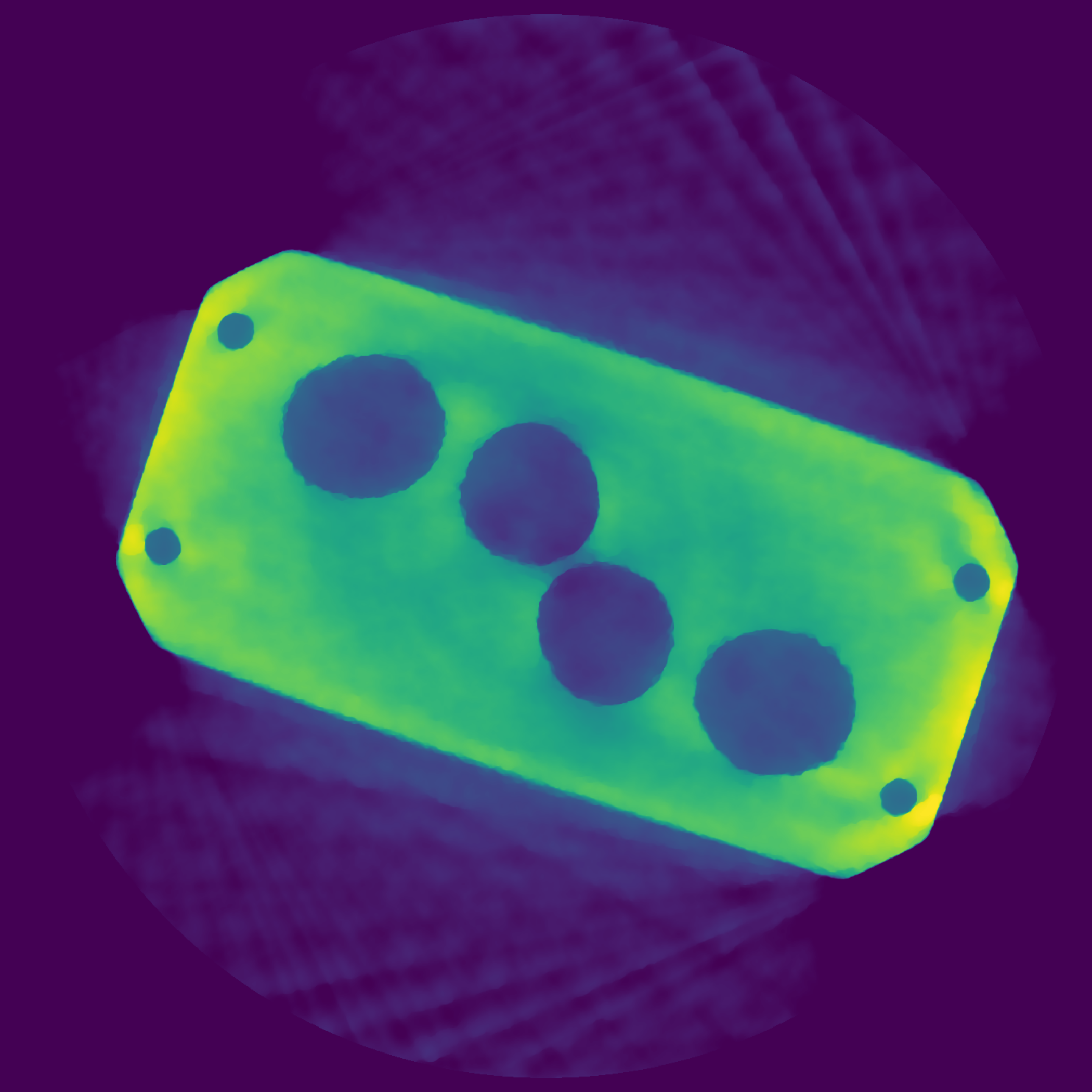}\\
    \includegraphics[width=0.25\linewidth]{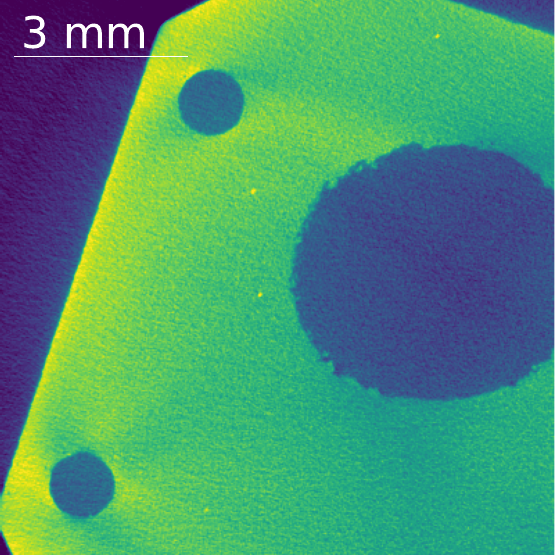}\hspace{-3pt}
    \includegraphics[width=0.25\linewidth]{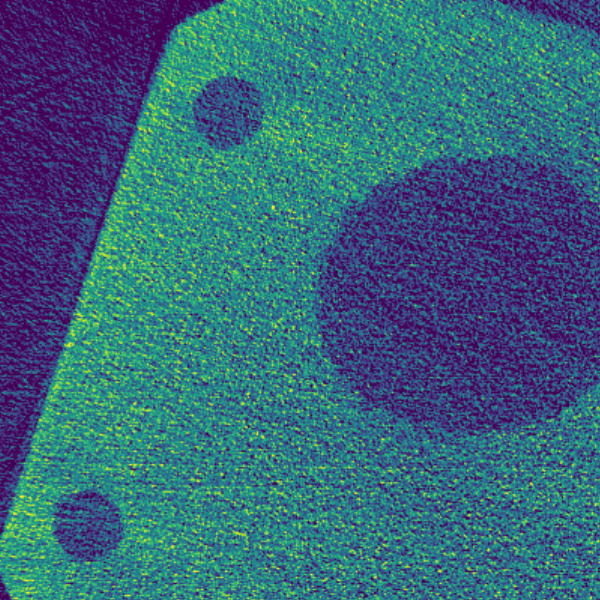}\hspace{-3pt}
    \includegraphics[width=0.25\linewidth]{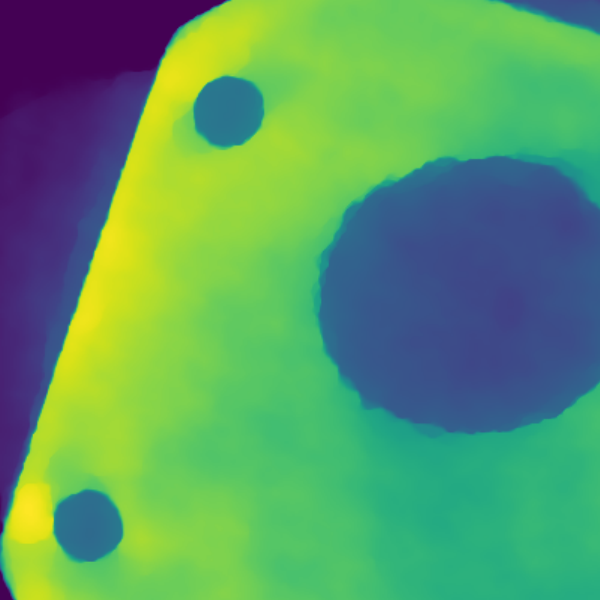}
    \caption{Industrial CT experiment. From left to right: ground-truth (FDK with 3142 projections), FDK with 60 projections, FISTA-aCV with 60 projections. Top: $1400^2$ pixels. Bottom: $400^2$ pixels.}
    \label{fig:recon}
\end{figure}

Figure \ref{fig:recon} shows the results of the ground truth reconstruction, the few-view FDK reconstruction and the FISTA-CV reconstruction with $\lambda$ estimated with FISTA-aCV. Finally, we do not claim that there is no other algorithm competitive with FISTA-CV, e.g., ADMM, PDHG, PD3O \cite{condat_review2023}, but these algorithms do not achieve the optimal rate of FISTA for (\ref{eq.optimizationgenearal}) or exhibit the low per-iteration complexity of CV for (\ref{eq.denoising}).

\section{Conclusions}

We have presented FISTA-aCV, a memory efficient algorithm for the regularization hyperparameter learning related to non-smooth variational problems. It is based on the algorithms FISTA and Condat-Vu within a bilevel learning problem addressed with AD to compute the derivatives requiring reduced memory storage. Then, an accelerated first-order algorithm is obtained to reach a solution. This work was motivated by few-view industrial CT and the methodology was tested with experimental tomographic data where both the parameter learning and the resulting image reconstruction were addressed.   

\subsection*{Acknowledgements}
\noindent The work of Patricio Guerrero was supported by a FWO Research Project under Grant G0C2423N

\bibliographystyle{abbrv} 
\bibliography{main}
\end{document}